\documentclass[10pt]{article}
\usepackage{latexsym}

\author{F. Boniver\footnote{Aspirant FNRS}~ and P. B. A. Lecomte}
\title{A remark about the Lie algebra of infinitesimal conformal
transformations of the Euclidian space}
\date{}
\advance\topmargin by -\headheight
\advance\topmargin by -\headsep

\def\R{{\rm I\!R}}
\def\C{{\rm C\mkern-9.5mu I\mkern3mu}}

\newtheorem {Def}{Definition}   

\newtheorem {Lem}[Def]{Lemma}
\newtheorem {The}[Def]{Theorem}

\newtheorem {Pro}[Def]{Proposition}

\newenvironment{proof}[1][]{\par\addvspace{1em}\par\textsl{Proof{#1}. }}%
{\par\addvspace{1em}}
\newenvironment{lemma}{\par
\addvspace{1em}
\begin{Lem}
}{%
\end{Lem}
\par
\addvspace{1em}
}

\newenvironment{theorem}{\par
\addvspace{1em}
\begin{The}
}
{%
\end{The}
\par
\addvspace{1em}
}

\newcommand{\I}{{\rm 1\mkern-4.5mu I\mkern3mu}}

\def\Lie#1{L_{#1}}              % d

\newcommand{\alg}[2]{\mbox{$\mathit{{#1}}({#2},\R)$}} % 

\newcommand{\cro}[2]{[{#1},{#2}]}     % c

\newcommand{\algeng}[1]{\mbox{$>{#1}<$}} % s

\newcommand{\dx}[1]{\mathit{dx}^{{#1}}}

\newcommand{\glnr}{\mathit{gl}(n,\R)}
\newcommand{\Rnd}{\R^{n*}}
\newcommand{\Rn}{\R^n}

\newcommand{\lalg}{\mathit{so}(p+1,q+1)}

\newcommand{\Vect}{\mathrm{Vect}}
\newcommand{\Vectpoly}[1][*]{\mathrm{Vect}_{#1}}

\newcommand{\der}[1]{\partial_{#1}}

\newcommand{\conf}[1]{\mathit{Conf}({#1})}
\newcommand{\confpoly}[1]{\mathit{Conf}_*({#1})}
\begin{document}
\maketitle
%\baselinestretch
\begin{abstract}
Infinitesimal conformal transformations of $\Rn$ are always polynomial and finitely generated when $n>2$.  Here we
prove that  the Lie algebra of infinitesimal conformal polynomial transformations over $\Rn$, $n\geq 2$, is maximal in the 
Lie algebra of polynomial vector fields.  When $n$ is greater than $2$ and $p,q$ are such that $p+q=n$, this implies the
maximality  of an embedding of
$\alg{so}{p+1,q+1}$ into polynomial vector fields that was revisited in recent works about equivariant quantizations. 
It also refines a
similar but weaker theorem by V. I. Ogievetsky.
\end{abstract}
\vspace{1in}
AMS classification numbers : 17B66, 53A30
\clearpage
\section{Introduction}
Interest has been shown in recent works about equivariant quantizations (see e.g. \cite{do, lo}) for 
some particular  Lie 
subalgebras of vector fields over the $n$-dimensional Euclidian space $\Rn$.
These are embeddings of $\alg{sl}{n+1}$ and
 $\alg{so}{p+1,q+1}$, $(p+q=n)$, into the Lie algebra of polynomial vector fields over $\Rn$.
% which we will denote by $\Vectpoly[*]$. 

 Because of its interpretation as the infinitesimal counterpart to the action 
of the group $\alg{SL}{n+1}$ on the $n$-dimensional real projective space, the embedding of $\alg{sl}{n+1}$
 into these fields is called the
\emph{projective Lie algebra}. 
It is quite a well-known fact that it is maximal among polynomial vector fields; see \cite{lo} for a proof.

In this paper, we focus our interest on the other subalgebras, related to infinitesimal conformal transformations. 

We will in particular generalize and refine a theorem by V. I. Ogievetsky: in \cite{og},  a
mix of projective and conformal vector fields is said to generate the Lie algebra of polynomial vector fields and an explicit
proof  given in dimension
$4$.

When dimension $n$ is greater than 2, infinitesimal
conformal transformations constitute a finite-dimensional Lie algebra, made up of polynomial vector
fields.  This is the considered embedding of $\alg{so}{p+1,q+1}$ into vector fields.   When $n=2$, this embedding is only
a finite-dimensional subalgebra of the (infinite-dimensional) Lie algebra of conformal infinitesimal transformations.

In both cases, we prove that the subalgebra of such polynomial conformal transformations is maximal among polynomial 
vector fields.  For the sake of completeness, we examine the special position of the introduced finite-dimensional 
subalgebras
in dimension
$2$.

We will denote by
$\Vect(\Rn)$ the Lie algebra of vector fields over
$\Rn$, and respectively by
$\Vectpoly(\Rn)$,
$\Vectpoly[\leq i](\Rn)$ and $\Vectpoly[i](\Rn)$ the Lie algebra of polynomial vector fields over $\Rn$, the space of
polynomial  fields of degree not greater than $i$ and the space of homogeneous  fields of degree $i$.

We will always assume that dimension $n$ is greater than $1$.

\section{The  algebra $\conf{p,q}$} 
Let us denote by $\conf{p,q}$ the Lie algebra of  vector fields over $\Rn$, $\mbox{$n=p+q$}$,  conformal with respect to the
metric
\[
g=\sum_{i=1}^n a_i (\dx{i})^2,
\]
where $a_1=\ldots=a_p=1$ and $a_{p+1}=\ldots=a_n=-1$.

These are the  fields $X$ which satisfy
\[
\Lie{X}g=\alpha_X g
\] for some smooth function $\alpha_X$.  Denote by $\der{i}$ both the partial derivative along the $i$-th axis and
the $i$-th natural basis vector of $\Rn$. It is equivalent for the components of
$X=\sum_i X^i\der{i}$ to satisfy
\begin{equation}\label{l1}
	\left\{
		\begin{array}{lll}
			\der{i}(a_j X^j)+\der{j}(a_i X^i)&=&0 \\
			\der{i}(a_i X^i)-\der{j}(a_j X^j)&=&0	
		\end{array}
	\right.
\end{equation}
when $i\neq j$.

For $h\in\Rn$, $A\in\glnr$ and $\alpha\in\Rnd$, define
\[
\begin{array}{lll}
h^* &=& -\sum_i h^i\der{i} \\
A^* &=& -\sum_{i,j} A^i_j x^j\der{i} \\
\alpha^* &=& \alpha(x)\sum_i x^i\der{i}-\frac{1}{2}(\sum_i a_i (x^i)^2)\alpha^\sharp
\end{array}
\]
where $\alpha^\sharp=\sum_i a_i \alpha_i \der{i}$.

The space 
\[
(\Rn\oplus\alg{so}{p,q}\oplus\R\I\oplus\Rnd)^*
\]
is a Lie subalgebra of $\Vect(\Rn)$,
isomorphic to $\alg{so}{p+1,q+1}$.  We will denote it by $\lalg$.

If $n>2$, it is known that $\conf{p,q}=\lalg$.  This follows for instance from \cite{ta}.

When $p=2$, $q=0$, condition $(\ref{l1})$  precisely means $X^1+i X^2$ is holomorphic in $\C$: $\conf{2,0}$ is 
then isomorphic to the Lie subalgebra 
\[\{f \frac{\mathrm{d}}{\mathrm{d}z}:f \mbox{ is holomorphic in $\C$} \}
\] of $\Vect(\C)$.  
Through this isomorphism, $\mathit{so}(3,1)$ is mapped onto $\Vectpoly[\leq 2](\C)$ and is isomorphic to $\mathit{sl}(2,\C)$
considered as a real Lie algebra.

When $p=q=1$, the classical change of coordinates
\[
\left\{
\begin{array}{lll}
x^1+x^2&=&2u^1\\
x^1-x^2&=&2u^2\\
\end{array}
\right.
\]
transforms $\conf{1,1}$ into the space of smooth vector fields
\[
U^1(x^1)\der{1}+U^2(x^2)\der{2}.
\]
In other words, $\conf{1,1}$ is isomorphic to $\Vect(\R)\times\Vect(\R)$.
This time, $\mathit{so}(2,2)$ is mapped onto $\Vect_{\leq 2}(\R)^2$.

The particular form of condition $(\ref{l1})$ implies the following result, which turns out to be immediate but useful.
\begin{lemma}\label{lemma1}
If $\der{1}X,\ldots,\der{n} X\in\conf{p,q}$ then $X\in\conf{p,q}+\Vect_1(\Rn)$.
\end{lemma}
\begin{proof}
The first order derivatives of $X$ are conformal if and only if left hand sides of $(\ref{l1})$ are constant.
Subtracting a linear vector field from $X$, one can force them to vanish and thus $X$ to be conformal.
\end{proof}

\section{Maximality of conformal algebras of polynomial vector fields}
Denote by $\confpoly{p,q}$ the Lie subalgebra of $\conf{p,q}$ made up of polynomial vector fields.
Here is the announced result.

\begin{theorem}\label{thmmax}
$\confpoly{p,q}$ is maximal in $\Vectpoly(\Rn)$.
\end{theorem}
The word maximal is to be taken in its usual algebraic sense : the
only  algebra of polynomial vector fields larger than
$\confpoly{p,q}$ is $\Vectpoly(\Rn)$.

In order to prove the theorem, it suffices to show that any larger subalgebra than $\confpoly{p,q}$ contains every constant,
linear or quadratic vector field, as implied by the following straightforward lemma. 
\begin{lemma}\label{lmv2}
The smallest Lie subalgebra containing $\Vectpoly[\leq 2](\Rn)$ is  $\Vectpoly[*](\Rn)$.
\end{lemma}
It is of course not true when $n=1$, as $\Vectpoly[\leq 2](\R)$ is a subalgebra of $\Vectpoly(\R)$.

\begin{proof}[ of theorem \ref{thmmax}]
Let $X$ be a polynomial vector field not in $\confpoly{p,q}$.
We may suppose that
$X\in\Vectpoly[1](\Rn)\setminus (\alg{so}{p,q}\oplus\R\I)^*$. Indeed, if it is not
the case, repeatedly applying lemma~\ref{lemma1}, we replace
$X$ by some
$\der{i}X=\cro{\der{i}}{X}\notin \confpoly{p,q}$ 
as long as  possible and then subtract from the last built field its homogeneous parts of degree $\neq 1$.

Besides, as a module of $\alg{so}{p,q}$, $\glnr$ is split into three irreducible components:
\[
\glnr=\R\I\oplus\alg{so}{p,q}\oplus\alg{so}{p,q}^+,
\]
where $\alg{so}{p,q}^+$ denotes the space of traceless self-conjugate matrices with respect to the metric.  The linear
map
\[
A\in\glnr\mapsto A^*\in\Vect_1(\Rn)
\]
being an isomorphism of Lie algebras, it follows that the iteration of brackets of fields of $\alg{so}{p,q}^*$ with $X$
allows to generate  every linear vector field, i.e.
\[
\Vectpoly[1](\Rn)\subset\algeng{\confpoly{p,q}\cup\{X\}},
\]
if $>S<$ denotes the smallest Lie algebra containing a set $S$ of vector fields.

We still need to show that the latter algebra contains every homogeneous quadratic vector field.  

It is known that the action $\cro{A^*}{\cdot}$ of $A\in\glnr$ endows $\Vectpoly[2]$ with a structure of
$\glnr$-module for which  the subspace of divergence-free vector fields is  irreducible. 
Now, writing a quadratic vector field
\[
X=\frac{1}{n}\alpha^*+(X-\frac{1}{n}\alpha^*)
\]
with $\alpha\in\Rnd$ such that $\alpha(x)=\mathrm{div} X$, one easily sees that
 \[
  \Vectpoly[2](\Rn)=(\Vectpoly[2](\Rn)\cap\confpoly{p,q})\oplus\ker{\mathrm{div}}.
  \] 
Therefore, to generate $\Vectpoly[2]$, it suffices to find $Y$, $Z$ among the fields generated so far
such that $\cro{Y}{Z}\neq 0$ and $\mathrm{div}\cro{Y}{Z}=0$.
The fields $Y=-x^1\der{1}$ and $Z=(\dx{2})^*$ do the job.

Hence the result.
\end{proof}

As a consequence, $\lalg$ is not maximal in $\Vectpoly(\Rn)$ only if $n=2$.  It is easy to check that
$\mathit{so}(3,1)$ is maximal in the Lie subalgebra
\[
\{X\in\Vectpoly(\R^2) : X^1+i X^2 \mbox{ is a polynomial of } z=x+i y\}
\]
and that $\mathit{so}(2,2)$ is maximal in two subalgebras isomorphic to 
\[
\Vectpoly(\R)\times\mathit{sl}(2,\R),
\] in turn
maximal in a copy of $\Vectpoly(\R)\times\Vectpoly(\R)$.
%%%%%%%%%% BIBLIOGRAPHY %%%%%%%%%%%%%%

\vspace{.5in}
\noindent
Institut de math\'ematique, B37\\
Grande Traverse, 12\\
B-4000 Sart Tilman (Li\`ege)\\
Belgium\\
mailto:f.boniver@ulg.ac.be\\
mailto:plecomte@ulg.ac.be

\end{document}